\documentclass[11pt]{article}
\usepackage{amsmath}
\usepackage{amscd}
\usepackage{epsfig}

\usepackage[latin1]{inputenc}

\input{amssym}

\newcommand{\color}[6]{}
\newcommand\QQ{\hbox{I\kern-.53em\hbox{Q}}}
\newcommand\qed{\hfill$\sqcap\kern-8.0pt\hbox{$\sqcup$}$}
\newcommand\NN{\hbox{I\kern-.2em\hbox{N}}}

\newcommand\RR{\hbox{I\kern-.2em\hbox{R}}}
\newcommand\sRR{{\sl \hbox{I\kern-.2em\hbox{R}}}}

\newcommand{\pp}{{\bf P}^1}

\newcommand\ZZ{{{\rm Z}\kern-.28em{\rm Z}}}
\newcommand\proof{\noindent{\em{Proof}.\ }}

\newcommand{\Bc}{T_{\textrm{\scriptsize bif}}}

\newtheorem{theo}{Theorem}[section]
\newtheorem{prop}[theo]{Proposition}
\newtheorem{lem}[theo]{Lemma}
\newtheorem{rem}[theo]{Remark}
\newtheorem{cor}[theo]{Corollary}

\newcommand{\la}{\lambda}

\numberwithin{equation}{section}

\begin{document}

\date{}
\title{Distribution of polynomials with cycles of given multiplier}

\vskip0.5cm

\author{Giovanni Bassanelli and
 Fran\c{c}ois Berteloot
}

\vskip0.5cm

\maketitle

\begin{abstract}
In the space of degree $d$ polynomials, the hypersurfaces
defined by the existence of a cycle of period $n$ and multiplier $e^{i\theta}$ are known to be contained in the bifurcation locus.
We prove that these hypersurfaces
equidistribute the bifurcation current. In particular, degree $d$ polynomials having a cycle of multiplier  $e^{i\theta}$ are dense in
the bifurcation locus; this is a new result even when $d=2$.
\end{abstract}

{\footnotesize 

Fran\c cois Berteloot, Universit\' e Paul Sabatier MIG. 
Institut de Math\'ematiques de Toulouse. 
31062 Toulouse
Cedex 9, France.
{\em Email: berteloo@picard.ups-tlse.fr}}

\medskip
\noindent Math Subject Class: 37F45; 37F10

\section{Introduction}

Within the space ${\cal P}_d$ of degree $d$ polynomials, the sets $Per_n(w)$ of those having a cycle of exact period $n$ and multiplier $w$ turn
out to be hypersurfaces. One knows, since the fundamental work of Ma\~n\'e, Sad and Sullivan \cite{MSS}, that the closure of the union of the hypersurfaces $Per_n(e^{i\theta})$
coincides with the bifurcation locus of ${\cal P}_d$, that is the set of polynomials whose dynamics drastically changes under small perturbation.
However, nothing is known about the asymptotic behaviour of $Per_n(e^{i\theta})$ as the period $n$ grows. We aim,
in the present article, to give a precise answer to this problem. \\

In the first part of the paper we actually deal with holomorphic families of rational maps $\big(f_{\la}\big)_{\la \in M}$. Our approach exploits the properties of the
Lyapunov function $L$ which assigns to each $\la\in M$ the Lyapunov exponent $L(\la)$ of $f_{\la}$ with respect to its maximal entropy measure.
The key is an approximation property for the Lyapunov exponent (see theorem \ref{theoapprox}), which naturally relates $L(\la)$ and $L_n(\la,w):=d^{-n} \ln \vert p_n(\la,w)\vert$ where the functions $p_n(\cdot,w)$  canonically define the hypersurfaces $Per_n(w)$
 (see theorem \ref{theopoly}). \\
We may thus compare the limits of $d^{-n}[Per_n(w)] :=dd^c\; L_n(\cdot,w)$ with the bifurcation current $T_{bif}$ since, as it has been shown by DeMarco (\cite{DeM2}), 
$T_{bif}$ coincides with the current $dd^c\;L$ (see also \cite{BB}). In this spirit, the theorem \ref{theoappL} gathers the results obtained by combining simple dynamical
properties with potential-theoretic ones.
Its first two statements are quite easily obtained and 
 were actually essentially given in our paper \cite{BB2}:

$$d^{-n}\;[Per_n(w)]\to T_{bif},\;\textrm{when}\;
\vert w\vert <1$$
$$\frac{d^{-n}}{2\pi}\int_0^{2\pi} [Per_n(re^{i\theta})]\;d\theta \to T_{bif},\;\textrm{when}\;
r\ge 0.$$

When $\vert w\vert \ge 1$ however, the convergence of $d^{-n}\;[Per_n(w)]$ is absolutely not clear. This is why we introduce some functions $L_n^+(\la,w)$
which coincide with $L_n(\la,w)$ when $\la$ is hyperbolic and do converge to $L$ in $L^1_{loc}(M)$ for any fixed $w\in {\bf C}$.
This leads to the following general equidistribution result on $M\times{\bf C}$
$$d^{-n}\;dd^c_{(\la,w)} \ln \vert p_n(\la,w)\vert\to T_{bif}.$$

To go further we assume that
 the hyperbolic parameters are sufficently well distributed in the parameter space and exploit the fact that, under such assumptions, the convergence of $L_n(\cdot,w)$ to $L$
in $L^1_{loc}$ for $\vert w\vert =1$ may be deduced from that of $L^+_n(\cdot,w)$. This leads to some results  which are then used for studying
polynomial families (see, in particular, the proposition \ref{propmodel}). Let us stress that if  hyperbolic parameters are dense in $M$
then $d^{-n}[Per_n (w)]\to T_{bif}$ for any $w\in {\bf C}$. By the work of Przytycki, Rivera-Letelier and Smirnov \cite{PRS}, this actually also occurs when Collet-Eckman parameters
are dense. \\

In the second part of the paper, we restrict ourselves to polynomials and prove the following:

\begin{theo} \label{theoequipol} Let $d\ge 2$ and $\{P_{c,a}\}_{ (c,a)\in{\bf C}^{d-1}}$ be the holomorphic family of degree $d$ polynomials parametrized 
 by defining $P_{c,a}$ as the polynomial of degree $d$ whose critical points are $(0,c_1,\cdot\cdot\cdot,c_{d-2})$
and such that $P_{c,a}(0)=a^d$. Let $T_{bif}$ be the bifurcation current of this family. Then $\lim_n d^{-n}[Per_n(w)]=T_{bif}$ for any $w$ such that $\vert w\vert\le 1$.\\
\end{theo}

Roughghly speaking, the proof of theorem \ref{theoequipol} consists in showing that $L$ is the only $L^1_{loc}$ limit value of $L_n(\cdot,e^{i\theta})$ by controling the convergence
on some slices where we may use our previous results (this is sketched at the beginning of subsection \ref{proof} in the case $d=3$). The existence of such slices, which are chosen for the good repartition of hyperbolic parameters, is intimately related to the behaviour at infinity of the bifurcation locus. More precisely, 
 we need to control, in a projective compactification of the parameter space ${\bf C}^{d-1}$, 
how the sets of parameters
$(c,a)$ for which a given critical point has a bounded orbit cluster at infinity. A precise description of these cluster sets (see theorem \ref{controlinfty})
actually follows from the work 
of Branner and Hubbard on the compactness of the connectedness locus.\\ 

To end this introduction, we would like to mention the works \cite{DF} of Dujardin and Favre and \cite{Du} of Dujardin where a different approach of the bifurcation current
is developed and, in particular,  distribution results for the Misiurewicz parameters are proved.

\section{Some tools}

\subsection{Hypersurfaces $Per_n(w)$}\label{sspn} 

For any holomorphic family of rational maps,
the following result precisely describes the set of maps
having a cycle of given period and multiplier.

\begin{theo}\label{theopoly}
Let $f:M\times\pp\to\pp$ be a holomorphic family of degree $d\ge 2$ rational maps. Then for every integer $n\in\NN^*$ there exists a
holomorphic function $p_n$ on $M\times {\bf C} $ which is polynomial on ${\bf C}$ and such that:
\begin{itemize}
\item[1-] for any $w \in {\bf C} \setminus \{1\}$, the function
$p_n(\la,w)$ vanishes if and only if $f_{\la}$ has a cycle of exact period $n$ and multiplipler $w$
\item[2-] $p_n(\la,1)=0$ if and only if $f_{\la}$ has a cycle of exact period $n$ and multiplier $1$ or a cycle of exact period $m$ whose multiplier is a primitive $r^{th}$ root of unity
with $r\ge 2$ and $n=mr$
\item[3-] for every $\la\in M$, the degree $N_d(n)$ of $p_n(\la,\cdot)$ satisfies $ d^{-n} N_d(n) \sim\frac{1}{n}$.
\end{itemize}
\end{theo}

This leads to the following definitions: for any integer $n$ and any $w\in {\bf C}$ the subset $Per_n(w)$ of $M$
is the hypersurface given by 
$$Per_n(w):=\{\la\in M/\;p_n(\la,w)=0\}$$
and, taking into account the possible multiplicities, we consider the following integration currents:
$$[Per_n(w)]:=dd^c_{\la}\;\ln\vert p_n(\la,w)\vert.$$

Let us briefly recall the construction of the functions $p_n$. For more details we refer to the paper of Milnor \cite{Mi} or the fourth chapter of the book
of Silverman \cite{Silbook}.\\

One first constructs the dynatomic polynomials $\Phi_{\varphi,n}^*$ associated to a rational map $\varphi$ of degree $d\ge 2$.
Let us denote by $F^n=(F^n_1,F^n_2)$ the iterates of some lift $F$ of $\varphi$ to ${\bf C}^2$ and define homogeneous polynomials $\Phi_{\varphi,n}$
on ${\bf C}^2$ by setting:

$$\Phi_{\varphi,n}(X,Y):=YF_1^n(X,Y)-XF_2^n(X,Y).$$

The divisor $Div(\Phi_{\varphi,n})$ induced by $\Phi_{\varphi,n}(X,Y)$ on $\pp$ is precisely the set of periodic points of $\varphi$ with exact period dividing $n$.
Denoting $\mu$ the classical M\"{o}bius function, one then sets 

$$\Phi_{\varphi,n}^*(X,Y):=\Pi_{k\vert n}\big(\Phi_{\varphi,n}(X,Y)\big)^{\mu(\frac{n}{k})}.$$

 Using the fact that the sum $\sum_{k\vert n}\mu(\frac{k}{n})$ vanishes if $n>1$ and is equal to $1$ if $n=1$, one may show that 
$\Phi_{\varphi,n}^*$ is a polynomial whose degree $\nu_d(n)$ depends only on 
$n$ and $d$. The divisor $Div(\Phi_{\varphi,n}^*)$ induced by $\Phi_{\varphi,n}^*(X,Y)$ on $\pp$ clearly contains the periodic points
of $\varphi$ with exact period equal to $n$. The other points contained in $Div(\Phi_{\varphi,n}^*)$ are precisely the periodic points of $\varphi$ whose
exact period $m$ divides $n$ ($m=nr,r\ge2$) and whose multiplier is a primitive $r$-th root of unity (see \cite{Silbook} theorem 4.5 page 151).\\ 
If $z\in Div(\Phi_{\varphi,n}^*)$ has exact period $m$ with $n=mr$,  we shall denote by $w_n(z)$ the $r$-th power of the multiplier of $z$ (that is 
$(\varphi^n)'(z)$ in good coordinates).
One sees in particular that the following fact occurs: 
a point $z$ is periodic of exact period $n$ and $w_n(z)\ne 1$ if and only if
$z\in Div(\Phi_{\varphi,n}^*)$ and $w_n(z)\ne 1$.

Let us now consider the sets

$$\Lambda_n^*(\varphi):=\{w_n(z);\;z\in Div(\Phi_{\varphi,n}^*)\}$$

where the points in $Div(\Phi_{\varphi,n}^*)$ are counted with multiplicity and let us denote by $\sigma_i ^{*(n)}(\varphi)$, $1\le i\le \nu_d(n)$,  the associated symmetric functions.
We define the polynomials $p_n(\varphi,w)$  by

$$\big(p_n(\varphi,w)\big)^n:=\Pi_{i=0}^ {\nu_d(n)} \sigma_i ^{*(n)}(\varphi) (-w)^ {\nu_d(n)-i}$$

and therefore $p_n(\varphi,w)=0$ if and only if $w\in\Lambda_n^*(\varphi)$. The properties of $p_n$ follow easily from this construction.
The degree $N_d(n)$ of $p_n(\la,\cdot)$ is equal to $\frac{1}{n}\nu_d(n)=\frac{1}{n}\sum_{k\vert n}\mu(\frac{n}{k})d^k$. In particular $ d^{-n} N_d(n) \sim \frac{1}{n}$.
\\

 \subsection{Lyapunov exponent and bifurcation current}\label{sstbif}
 Every rational map of degree $d\ge 2$ on the Riemann sphere admits a maximal entropy measure
$\mu_f$. 
The Lyapunov exponent of $f$ with respect to the measure $\mu_f$ is given by $L(f)=\int_{{\bf P}^1} \log \vert f'\vert\mu_{f}$
(see \cite{DS} for a general exposition in any dimension).\\

When $f:M\times {\bf P}^1\to {\bf P}^1$ is a holomorphic family of degree $d$ rational maps, the Lyapunov function $L$ on the parameter space $M$
is defined by:
$$L(\la)=\int_{{\bf P}^1} \log \vert f_{\la}'\vert\mu_{\la}$$
where $\mu_{\la}$ is the maximal entropy measure of $f_{\la}$. The function $L$ is $p.s.h$ on $M$ and the bifurcation current $\Bc$ of the family is a closed, positive $(1,1)$-current on $M$ which may be defined by
$$\Bc:=dd^c L(\la).$$
As it has been shown by DeMarco \cite{DeM2}, the support of $\Bc$ concides with the bifurcation locus of the family  in the sense of 
Ma\~n\'e-Sad-Sullivan (see also \cite{BB}, Theorem 5.2).\\

Let us recall that Ma\~n\'e, Sad and Sullivan have shown that the complement of the bifurcation locus is a dense open subset of the parameter space $M$
whose connected components are the so-called \emph{stable components}. They have also shown that any neutral cycle is persistent on the stable components, in the language of theorem 
\ref{theopoly}, this property may be expressed as follows:

\begin{rem}\label{compstab}  For $\vert w_0\vert =1$, a function
$p_n(\la,w_0)$ either does not vanish on any stable component or vanishes identically on $M$.
\end{rem}

In our study we shall combine classical potential-theoretic methods with the following 
dynamical property ( see \cite{BDM} where this has been proved for endomorphisms of ${\bf P}^k$).

\begin{theo}\label{theoapprox}
 Let $f:\pp\to\pp$ be a degree $d\ge 2$ rational map,  $\mu$ its maximal entropy measure and $L$ the Lyapunov exponent of $f$ with respect to $\mu$.
Then:
$$L=\lim_n \frac{d^{-n}}{n} \sum_{p\in R^{*}_n} \ln \vert (f^n)'(p)\vert$$
where $R^{*}_n:=\{p\in \pp \;/\; \;p\; \textrm{has exact period}\; n \;  \textrm{and}\;\vert (f^{n })' (p)\vert > 1 \}$.
\end{theo}

The continuity of the Lyapunov function will also play a crucial role. This was proved by Ma\~n\'e \cite{Ma} but a simple argument based on DeMarco's formula
shows that this function is actually H\"{o}lder continuous (\cite{BB} Corollary 3.4). 

\begin{theo}\label{theocontL}
 Let $f:M\times\pp\to\pp$ be a holomorphic family of degree $d\ge 2$ rational maps. Let $L(\la)$ be the Lyapunov exponent of $(\pp,f_{\la},\mu_{\la})$ where $\mu_{\la}$ is the maximal entropy measure of 
$f_{\la}$. Then the function $L(\la)$ is H\"{o}lder continuous on $M$.
\end{theo}

\subsection{Potential theoretic tools}

The results from potential theory which we shall use are classical.
The main one concerns compacity properties of subharmonic functions (the second statement is known as
Hartogs lemma):

\begin{theo}\label{compapsh}
Let $(\varphi_j)$ be a sequence of subharmonic functions  which is locally uniformly bounded from above on some domain $\Omega \subset \RR^n$. 
\begin{itemize}
\item[1-] If $(\varphi_j)$ does not converge
to $-\infty$ then a subsequence $(\varphi_{j_k})$ is converging in $L^1_{loc}(\Omega)$ to some subharmonic function $\varphi$. 
In particular, $(\varphi_j)$ is converging in $L^1_{loc}(\Omega)$ to some subharmonic function $\varphi$ if it is pointwise converging to $\varphi$.
\item[2-] If $(\varphi_j)$ is converging in $L^1_{loc}(\Omega)$ to some subharmonic function $\varphi$ then, 
for any compact $K\subset \Omega$ and any continuous function $u$ on $K$, one has 
$\limsup_k \;\sup_K (\varphi_{j_k}-u) \le \sup_K(\varphi - u).$
\end{itemize}
\end{theo}

We shall also need the following continuity principle:

\begin{theo}\label{ContPrin}
Let $\varphi$ be a subharmonic function on some Riemann surface $M$. If the restriction of $\varphi$ to the support of its Laplacian is continuous then $\varphi$
is continuous on $M$.
\end{theo}

\section{Distribution of $Per_n(w)$ in general families}

In this section we consider an arbitrary holomorphic family $f:M\times\pp\to\pp$ of degree $d\ge 2$ rational maps. 
We investigate the convergence of the currents $\frac{1}{d^n}[Per_n(w)]$ towards the bifurcation current $T_{bif}$ by considering the sequences of their potentials and, therefore, compare the Lyapunov function $L$ with the limits of $\frac{1}{d^n}\ln \vert p_n(\la,w)\vert$ where $p_n(\la,w)$ are the polynomials given
by theorem \ref{theopoly}. \\

This leads us to consider the following sequences of $p.s.h$ functions:
\begin{center} 
$\displaystyle L_n^r(\la):=\frac{d^{-n}}{2\pi}\int_0^{2\pi} \ln \vert p_n(\la,re^{i\theta})\vert\;d\theta$\\
\medskip
$\displaystyle L_n^+(\la,w):=d^{-n}\sum_{j=1}^{N_d(n)} \ln^+\vert w- w_{n,j} (\la)\vert$\\
\medskip
$\displaystyle L_n(\la,w):=d^{-n}\ln \vert p_n(\la,w)\vert$
\end{center}

where $p_n(\la,w)=:\Pi_{j=1}^{N_d(n)}(w-w_{n,j}(\la))$ are the polynomials associated to the family $f$
by the theorem \ref{theopoly}.\\

The pointwise convergence of $L_n(\la,w)$ to $L$ for $\vert w\vert <1$ is quite a straightforward consequence of theorem \ref{theoapprox} and immediately implies that
$d^{-n}[Per_n(w)]$ converges to $T_{bif}$ when $\vert w\vert <1$.
However, when $\vert w\vert \ge 1$ and $\la$ is a non-hyperbolic parameter, the control of $L_n(\la,w)=d^{-n}\sum\ln\vert w- w_{n,j}(\la)\vert$ is very delicate because $f_\la$ may have many cycles whose multipliers are close to $w$.
This is why we introduce the $p.s.h$ functions $ L_n^+$ which both coincide with $L_n$ on the hyperbolic components and are quite easily seen to converge nicely.
Our main result is the:

\begin{theo}\label{theoappL} Let $f:M\times\pp\to\pp$ be a holomorphic family of degree $d\ge 2$ rational maps. Let $L(\la)$ be the Lyapunov exponent of $(\pp,f_{\la},\mu_{\la})$ where $\mu_{\la}$ is the maximal entropy measure of 
$f_{\la}$. Let $(L_n)_n$, 
$(L^r_n)_n$ and 
$(L^+_n)_n$ be the sequences of $p.s.h$ functions defined as above.
Then:
\begin{itemize}
\item[1-] the sequence $L_n$ converges pointwise to $L$ on $M\times \Delta$ and, for any $w\in\Delta$, the sequence $L_n(\cdot,w)$ converges in $L_{loc}^1$ to $L$ on $M$.
\item[2-] The sequence $L^r_n$ converges pointwise and in $L_{loc}^1$ to $L$ on $M$ for $r\ge0$.
\item[3-] The sequence $L^+_n$ converges pointwise and in $L_{loc}^1$ to $L$ on $M\times{\bf C}$; for every $w\in {\bf C}$ the sequence $L^+_n(\cdot,w)$ converges in  $L_{loc}^1$ to $L$ on $M$. 
\item[4-] The sequence $L_n$ converges in $L_{loc}^1$ to $L$ on  $M\times{\bf C}$.
\end{itemize}
\end{theo}

Let us stress that the various convergence properties for $d^{-n}[Per_n(w)]$ given in the introduction follow immediately from the first, second and last statements
of the above theorem by taking $dd^{c}$.\\

\proof 
All the statements are local and therefore, taking charts, we may assume that $M={\bf C}^k$.
We write the polynomials $p_n$ as follows :$$p_n(\la,w)=:\Pi_{i=1}^{N_d(n)}\big(w-w_{n,j}(\la)\big).$$
Throughout the proof we shall use the fact that
$d^{-n}N_d(n)\sim\frac{1}{n}$ (see theorem \ref{theopoly}).
In particular, this implies that the sequences $L_n$ and $L_n^+$ are locally uniformly bounded from above.\\

\noindent$\bullet$ We first establish the {\it convergence of  $L_n(\la,w)$ when $\vert w\vert <1$.}
 According to theorem \ref{theopoly}, the set $\{w_{n,j}(\la)\; /\;w_{n,j}(\la)\ne 1\}$ coincides with the set of multipliers
of cycles of exact period $n$ (counted with multiplicity) from which the cycles of multiplier $1$ are deleted.  Using the notation
$R^{*}_n(\la):=\{p\in \pp \;/\; \;p\; \textrm{has exact period}\; n \;  \textrm{and}\;\vert (f_{\la}^{n })' (p)\vert > 1 \}$ we thus have

\begin{eqnarray}\label{A1}
\sum_{j=1}^{N_d(n)} \ln^+ \vert w_{n,j}(\la)\vert =\frac{1}{n} \sum_{p\in R^{*}_n(\la)} \ln \vert (f^n)'(p)\vert.
\end{eqnarray}
Since $f_\la$ has a finite number of non-repelling cycles (Fatou' theorem), one sees that there exists $n(\la) \in\NN$ such that
\begin{eqnarray}\label{A2}
n\ge n(\la) \Rightarrow \vert w_{n,j}(\la)\vert > 1,\;\textrm{for any}\;1\le j\le N_d(n).
\end{eqnarray}

By \ref{A1} and \ref{A2}, one gets
\begin{eqnarray*}
L_n(\la,0)=d^{-n}\sum_{j=1}^{N_d(n)} \ln \vert w_{n,j}(\la)\vert =d^{-n}\sum_{j=1}^{N_d(n)} \ln^+ \vert w_{n,j}(\la)\vert =\frac{d^{-n}}{n} \sum_{R^{*}_n(\la)} \ln \vert (f^n)'(p)\vert
\end{eqnarray*}

 for $n\ge n(\la)$ which, by theorem \ref{theoapprox}, yields:

\begin{eqnarray}\label{A3}
\lim_n L_n(\la,0)=L(\la),\; \forall \la \in M.
\end{eqnarray}

Let us now pick $w\in \Delta$. By \ref{A2},
 we have $L_n(\la,w)-L_n(\la,0)=d^{-n}\sum_j\ln \frac{\vert w_{n,j}(\la) - w\vert }{\vert w_{n,j}(\la) \vert}$ 
 and 
$\ln(1-\vert w\vert)\le\ln\frac{\vert w_{n,j}(\la) - w\vert }{\vert w_{n,j}(\la) \vert}\le  \ln(1+\vert w\vert)$ for $1\le j\le N_d(n)$ and $n\ge n(\la)$ .
We thus get 
\begin{eqnarray*}
d^{-n}N_d(n)\ln(1-\vert w\vert)\le\vert L_n(\la,w)-L_n(\la,0) \vert \le d^{-n}N_d(n) \ln(1+\vert w\vert)
\end{eqnarray*}
for $n\ge n(\la)$ and, using \ref{A3}, $\lim_n L_n(\la,w)=L(\la)$ for any $(\la,w) \in M\times {\Delta}$.\\
The $L_{loc}^1$ convergence of $L_n(\cdot,w)$ now follows immediately from theorem \ref{compapsh}.\\

\noindent$\bullet$ Let us show that the convergence of $L_n(\la,0)=L_n^0 $ implies the {\it convergence of $ L_n^r $} for any $r>0$.
We essentially will show that $\lim_n \vert L_n^r(\la) - L_n(\la,0)\vert =0$ by 
using the formula $\ln Max(\vert a \vert,r)=\frac{1}{2\pi}\int_0^{2\pi}\ln \vert a-re^{i\theta}\vert d\theta$.
Indeed, this formula yields

\begin{eqnarray*}
L_n^r(\la)=\frac{1}{2\pi d^n}\int_0^{2\pi} \ln \Pi_j\vert re^{i\theta}-w_{n,j}(\la)\vert d\theta=\\
d^{-n}\sum_j \ln Max(\vert w_{n,j}(\la)\vert,r).
\end{eqnarray*}

Since $\vert w_{n,j}(\la)\vert\ge 1$ for $n\ge n(\la)$ (see \ref{A2}), we deduce from the above identity that:

\begin{eqnarray*}
L_n^r(\la)=d^{-n}\sum_j \ln \vert w_{n,j}(\la)\vert +d^{-n}\sum_{1\le\vert w_{n,j}(\la)\vert < r} \ln \frac{r}{\vert w_{n,j}(\la)\vert }=\\
L_n(\la,0) + d^{-n}\sum_{1\le\vert w_{n,j}(\la)\vert < r}  \ln \frac{r}{\vert w_{n,j}(\la)\vert }
\end{eqnarray*}

and thus

\begin{eqnarray*}
0\le L_n^r(\la) -
L_n(\la,0)= d^{-n}\sum_{1\le\vert w_{n,j}(\la)\vert < r}  \ln \frac{r}{\vert w_{n,j}(\la)\vert }\le
d^{-n}N_d(n)\ln^+ r.
\end{eqnarray*}

By \ref{A3}, this implies that $L_n^r$ is pointwise converging to $L$. It also shows that
 $(L_n^r)_n$ is locally uniformly bounded from above which, by theorem \ref{compapsh}, implies
 that $(L_n^r)_n$ converges to $L$ in $L^1_{loc} (M)$.\\

\noindent$\bullet$ Let us now deal with the {\it convergence of $L_n^+$.} We will show that $L_n^+(\cdot,w)$ is pointwise converging to $L$ on $M$ for every $w\in{\bf C}$. As
$(L_n^+)_n$ is locally uniformly bounded, this implies the convergence of $L_n^+(\cdot,w)$ in $L_{loc}^1(M)$ (theorem \ref{compapsh}) and the convergence of $L_n^+$ in 
$L_{loc}^1(M\times{\bf C})$ then follows by Lebesgue theorem.\\
 We have to estimate $L^+_n(\la,w)-L_n(\la,0)=:\epsilon_n(\la,w)$ on $M$. Let us fix $\la \in M$, $w\in {\bf C}$ and pick $R>\vert w\vert $. We may asume that $n\ge n(\la)$ so that $\vert  w_{n,j}(\la)\vert\ge 1$
for all $1\le j\le N_d(n)$  (see \ref{A2}) and 
then decompose $\epsilon_n(\la,w)$ in the following way:

\begin{eqnarray*}
\epsilon_n(\la,w)=d^{-n}\sum_{1\le\vert w_{n,j}(\la)\vert<R+1} \ln^+\vert w_{n,j}(\la) - w\vert 
+ d^{-n}\sum_{\vert w_{n,j}(\la)\vert\ge R+1} \ln\frac{\vert w_{n,j}(\la) - w\vert }{\vert  w_{n,j}(\la)\vert }\\
-d^{-n}\sum_{1\le\vert w_{n,j}(\la)\vert<R+1} \ln \vert w_{n,j}(\la) \vert .
\end{eqnarray*}

We may write this decomposition as
$
\epsilon_n(\la,w)=:\epsilon_{n,1}(\la,w) + \epsilon_{n,2}(\la,w) - \epsilon_{n,1}(\la,0).
$
As $L_n(\la,0)$ is converging to $L$, 
we simply have to check that $\epsilon_{n,1}(\la,w)$ and $\epsilon_{n,2}(\la,w)$ tends to $0$ when $n$ tends to $\infty$.
One clearly has
$
0\le \epsilon_{n,1}(\la,w)\le d^{-n} N_d(n) \ln\big(2R+1\big)
$
and thus
$\lim_n  \epsilon_{n,1}(\la,w)=0$.
Similarly, 
$
\lim_n \epsilon_{n,2}(\la,w)=0
$ follows from the fact that,  
for $\vert w_{n,j}(\la)\vert >R+1>\vert w\vert +1$, one has:
\begin{eqnarray*}
\ln(1-\frac{R}{R+1})\le\ln\frac{\vert w_{n,j}(\la) \vert -R}{\vert w_{n,j}(\la) \vert}\le\ln\frac{\vert w_{n,j}(\la) - w\vert }{\vert w_{n,j}(\la) \vert}\le\ln\frac{\vert w_{n,j}(\la) \vert +R}{\vert w_{n,j}(\la) \vert}\le \ln(1+\frac{R}{R+1}).
\end{eqnarray*}
\\

\noindent$\bullet$ We are finally ready to prove the {\it $L_{loc}^1$ convergence of $(L_n)_n$.} 
As
 the functions $L_n$ are $p.s.h$ and the sequence $(L_n)_n$ is locally uniformly bounded from above, we shall again
use the compacity properties of $p.s.h$ functions given by theorem \ref{compapsh}. Since $L_n(\la,0)$ is converging to $L(\la)$, the sequence
$(L_n)_n$ is not converging to $-\infty$ and it therefore suffices to show that,  among $p.s.h$ functions on $M\times {\bf C}$, the function $L$ is the only possible
limit for  $(L_n)_n$ in $L^1_{loc} (M\times {\bf C})$.

Let $\varphi$ be a $p.s.h$ function on $M\times {\bf C}$ and $(L_{n_j})_j$ a subsequence of $(L_n)_n$  which converges to $\varphi$
in $L^1_{loc} (M\times {\bf C})$. Pick $(\la_0,w_0)\in M\times{\bf C }$. We have to prove that $\varphi(\la_0,w_0)= L(\la_0)$.\\

Let us first observe that $\varphi(\la_0,w_0)\le L(\la_0)$.
Take a ball $B_{\epsilon}$ of radius $\epsilon$ and centered at $(\la_0,w_0)\in M\times{\bf C }$, by the submean value property and the $L_{loc}^1$-
convergence
of $L_n^+$ we have:

\begin{eqnarray*}
\varphi(\la_0,w_0)\le \frac{1}{\vert B_{\epsilon}\vert}\int_{B_{\epsilon}}\varphi\; dm=\lim_j \frac{1}{\vert B_{\epsilon}\vert}\int_{B_{\epsilon}} L_{n_j}\;dm\\
\le \lim_j \frac{1}{\vert B_{\epsilon}\vert}\int_{B_{\epsilon}} L_{n_j}^+\;dm = \frac{1}{\vert B_{\epsilon}\vert}\int_{B_{\epsilon}} L\;dm
\end{eqnarray*}

and then, making $\epsilon\to 0$, one obtains $\varphi(\la_0,w_0)\le L(\la_0)$.\\
Let us now check that $\limsup_j L_{n_j}(\la_0,w_0e^{i\theta})=L(\la_0)$ for almost all $\theta\in [0,2\pi]$. Let $r_0:=\vert w_0\vert$. 
As we saw, $L_n^+$ is pointwise converging to $L$ and therefore:
$$\limsup_j L_{n_j}(\la_0,w_0e^{i\theta})\le \limsup_j L_{n_j}^+(\la_0,w_0e^{i\theta})=L(\la_0)$$

on the other hand, by pointwise convergence of $L_{n}^{r_0}$ to $L$ and Fatou's lemma we have:

\begin{eqnarray*}
L(\la_0)=\lim_nL_{n}^{r_0} (\la_0)=\limsup_j\frac{1}{2\pi}\int_0^{2\pi} L_{n_j}(\la_0,r_0 e^{i\theta}) d\theta \le\\
\frac{1}{2\pi}\int_0^{2\pi} \limsup_j L_{n_j}(\la_0,r_0 e^{i\theta}) d\theta. 
\end{eqnarray*}

and the desired property follows immediately.\\
To end the proof we argue by contradiction and assume that $\varphi(\la_0,w_0) < L(\la_0)$.
As $\varphi$ is upper semi-continuous and $L$ continuous (theorem \ref{theocontL} ), there exists a  neighbourhood $V_0$ of 
$(\la_0,w_0)$ and $\epsilon>0$ such that 
$$\varphi - L \le -\epsilon\;\;\textrm{on}\;V_0.$$
Pick a small ball $B_{\la_0}$ centered at $\la_0$ and a smal disc $\Delta_{w_0}$ centered at $w_0$ such that  
$B_0:=B_{\la_0}\times\Delta_{w_0}$ is relatively compact in $V_0$.  Then, according to Hartogs lemma (see theorem \ref{compapsh})
we have:

$$\limsup_j \big(Sup_{B_0} (L_{n_j}-L)\big)\le Sup_{B_0} (\varphi -L)\le -\epsilon.$$
This is impossible since, as we have seen before, we may find $(\la_0,r_0 e^{i\theta_0})\in B_0$ such that
$\limsup_j\big(L_{n_j}(\la_0,r_0e^{i\theta_0})-L(\la_0)\big)=0$.\qed\\

\begin{rem} Using standard techniques, one may deduce from the fourth assertion of theorem \ref{theoappL} that the set of multipliers $w$ for which
the bifurcation current $T_{bif}$ is not a limit of the sequence $d^{-n} [Per_n(w)]$ is contained in a polar subset of the complex plane.
\end{rem}

The fact that the functions $L_n^+$ and $L_n$ coincide on hyperbolic components would easily yield the convergence of 
$d^{-n}[Per_n(w)]$ towards $T_{bif}$ for any $w\in {\bf C}$ if the density of hyperbolic parameters in $M$ would be known.
Using the equivalence between \emph{uniform hyperbolicity on periodic orbits} and Collet-Eckmann conditions (see \cite{PRS}), one sees that the same conclusion occurs when some kind of
non-uniform hyperpolic parameters are dense.
Without knowing that, we are not able to prove this convergence when $\vert w\vert \ge 1$ but we can
overcome the difficulty when the hyperbolic parameters are sufficently nicely distributed. We establish now a few basic facts of this nature which will be used in our study of polynomial families
in the next section.\\

The following proposition summarizes some usefull remarks.

\begin{prop}\label{propelem} Let us make the same assumptions and adopt the same notations than in theorem \ref{theoappL}. Let $w_0\in {\bf C}$.
Then:
\begin{itemize}
\item[1-] any $p.s.h$ limit value of $L_n(\la,w_0)$ in $L^1_{loc}(M)$ is smaller than $L$.
\item[2-] If a subsequence $L_{n_k}(\la,w_0)$ converges pointwise to $L$ on the stable set  then it also converges to $L$ in $L^1_{loc}(M)$. 
\item[3-] Assume that $\vert w_0\vert=1$ and that the family has no persistent neutral cycle. 
If a subsequence $L_{n_k}(\la,w_0)$ converges to $\varphi$ in $L^1_{loc}(M)$ then $\varphi$ is pluriharmonic on any stable component $\Omega$ and the convergence is locally uniform on 
$\Omega$ . 
\end{itemize}
\end{prop}

\proof  1- Let us set $\varphi_n(\la):=L_n(\la,w_0)$ and assume that a subsequence $\varphi_{n_j}$ converges in $L^1_{loc}(M)$ to some $p.s.h$ function $\varphi$.
Since $\L^+_n(\la,w_0)$ converges to $L$ in $L^1_{loc}(M)$ and $\varphi_{n_j}(\la)\le \L^+_{n_j}(\la,w_0)$ we get
$\varphi(\la_0)\le \frac{1}{\vert B_{\epsilon}\vert}\int_{B_{\epsilon}}\varphi\; dm \le \frac{1}{\vert B_{\epsilon}\vert}\int_{B_{\epsilon}} L\;dm$ for any small ball $B_{\epsilon}$ centered at 
$\la_0$. The desired inequality then follows by making $\epsilon\to 0$.\\
2- Recall that the stable set is an open dense subset of $M$. Let $\varphi$ be any $p.s.h$ limit of $L_{n_k}(\la,w_0)$ in $L^1_{loc}(M)$, we have to show that $\varphi=L$.
By the first assertion $\varphi\le L$. As $\varphi=L$ on a dense subset, the semicontinuity of $\varphi$ and the continuity of $L$ (see theorem \ref{theocontL}) imply that $\varphi \ge L$.\\
3- Using the remark \ref{compstab} one sees that the functions $L_{n_k}(\la,w_0)$ are pluriharmonic on $\Omega$, this implies that
$\varphi$ itself is pluriharmonic on $\Omega$ and that $L_{n_k}(\la,w_0)$ is actually converging locally uniformly on $\Omega$ to $\varphi$.\qed\\

Let us precise how the density of hyperbolic parameters allows to strenghten the conclusion of theorem \ref{theoappL}.

\begin{prop}\label{prophyp} Let us make the same assumptions and adopt the same notations than in theorem \ref{theoappL}. Let $w_0\in {\bf C}$.
Then:
\begin{itemize}
\item[1-] for any hyperbolic component $\Omega\subset M$, the sequence $L_n(\la,w_0)$ converges locally uniformly to $L$ on $\Omega$.
\item[2-] If the hyperbolic parameters are dense in $M$ then the sequence $L_n(\la,w_0)$ converges to $L$ in $L^1_{loc}(M)$.
\end{itemize}
\end{prop}

\proof
1- If $\la$ is a hyperbolic parameter then $f_{\la}$ has only attracting or repelling cycles and is expansive on its Julia set. Thus, as $f_{\la}$ has at most a finite number of attracting cycles,
one sees that $\vert w_{n,j}(\la)\vert \ge \vert w_0\vert +1$ for all $1\le j\le N_d(n)$ and $n$ big enough. In other words $L_n(\la,w_0)=\L^+_n(\la,w_0)$ for $n$ big enough and therefore,
according to theorem \ref{theoappL}, $L_n(\la,w_0)$ converges to $L(\la)$. By theorem \ref{compapsh}, $L_n(\la,w_0)$ converges to $L$ in $L^1_{loc}(\Omega)$. The local uniform convergence then 
follows
from the previous proposition.\\
2- This follows immediately from an argument of unic limit value based on theorem \ref{compapsh} after combining the above assertion with the second one of proposition \ref{propelem}.\qed\\

In the remaining of the paper we will focus on the case $\vert w_0\vert=1$ and work with polynomial families. Slicing the parameter space in different ways, we will obtain one dimensional holomorphic families for which the problem is easier to handle. The following technical lemma covers the different situations which we shall consider.

\begin{lem}\label{lemmodel}  Let $M$ be a Riemann surface and $f:M\times\pp\to\pp$ be a holomorphic family of degree $d\ge 2$ rational maps.
Let $L(\la)$ and $L_n(\la,w)$ be the subharmonic functions defined in theorem \ref{theoappL}.
Let  $w_0\in {\bf C}$ with $\vert w_0\vert =1$ and $\varphi$ be a subharmonic limit value of $L_n(\la,w_0)$ in $L^1_{loc}(M)$ such that:
\begin{itemize}
\item[1-] the bifurcation locus is contained in the closure of  the set of parameters where $\varphi=L$
\item[2-] $\varphi=L$ on the stable component which are not relatively compact in $M$.
\end{itemize}
Then $\varphi=L$.
\end{lem}

\proof
We shall use several times the fact that the function $L$ is continuous (see theorem \ref{theocontL}).
Assume that $\varphi_{n_j}:=L_n(\cdot,w_0)$
converges to $\varphi$, then the holomorphic functions $p_{n_j}(\la,w_0)$ cannot vanish identically for $j$ big enough. According to the remark \ref{compstab} this implies that the functions $\varphi_{n_j}$
are harmonic on all stable components of $M$. This leads to the 
simple, but crucial, observation that $\varphi$ is harmonic on any stable component or, in other words, that the Laplacian $\Delta \varphi$ is supported in the bifurcation locus. \\
According to the first assertion of proposition \ref{propelem}, we have $\varphi\le L$. 
We may now see that $\varphi =L$ on the bifurcation locus. Indeed, if $\la_0$ belongs to the bifurcation locus then, by assumption, there exists a sequence $\la_k$ converging to $\la_0$ such that $\varphi(\la_k) =L(\la_k)$. Then, using the upper-semicontinuity of $\varphi$ and the continuity of $L$, we get $\varphi(\la_0) =\limsup_{\la\to\la_0} \varphi(\la) \ge \limsup \varphi(\la_k) =\lim L(\la_k) =L(\la_0)$.\\
By the first observation and the fact that $L$ is continuous, we see that $\varphi$ is continuous on the support of its Laplacian, this implies
that $\varphi$ is continuous on $M$ (see theorem \ref{ContPrin}). We may now prove that $\varphi\equiv L$. If this would not be the case, then $\varphi(\la_0)<L(\la_0)$
for some $\la_0\in M$. As $L$ and $\varphi$ coincide on the bifurcation locus and (by assumption) on non-relatively compact stable components, $\la_0$ would belong to some stable component $\Omega$
which is relatively compact in $M$. This contradicts the maximum principle since $(\varphi -L)$ is continuous on $\overline{ \Omega}$, harmonic on $\Omega$ and vanishes on $b\Omega$.\qed\\

Here is a typical application of the above lemma to the case where the hyperbolic parameters are well distributed in $M$, it covers the case of the quadratic polynomial family.

\begin{prop}\label{propmodel}  Let $M$ be a Riemann surface and $f:M\times\pp\to\pp$ be a holomorphic family of degree $d\ge 2$ rational maps which satisfies the two following conditions:
\begin{itemize}
\item[1-] the bifurcation locus is contained in the closure of  hyperbolic parameters
\item[2-] the set of non-hyperbolic parameters is compact in $M$.
\end{itemize}
Let $L(\la)$ and $L_n(\la,w)$ be the subharmonic functions defined in theorem \ref{theoappL}. Then, if $\vert w_0\vert =1$, the sequence $L_n(\la,w_0)$ converges to $L$ in $L_{loc}^1(M)$.
\end{prop}

\proof By the first assertion of proposition \ref{prophyp},  the sequence $L_n(\la,w_0)$ does not converge to $-\infty$.
According to the theorem \ref{compapsh} it thus suffices to show that 
any subharmonic limit value $\varphi$ of $L_n(\la,w_0)$ in $L^1_{loc}(M)$ coincides with $L$.
This follows immediately from lemma \ref{lemmodel} since,
by the first assertion of proposition \ref{prophyp} again, $\varphi=L$ on the non relatively compact stable components.\qed

\section{Distribution of $Per_n(w)$ in polynomial families}

\subsection{The space of degree $d$ polynomials}

Let ${\cal P}_d$ be the space of polynomials of degree $d\ge 2$ with $d-1$ marked critical points up to conjugacy by affine transformations. Although this space has a natural structure of affine variety of dimension $d-1$,
we shall actually work with a specific parametrization of ${\cal P}_d$ which has been introduced by Dujardin and Favre in \cite{DF}. We refer to their paper and to the seminal paper
\cite{BH} of Branner and Hubbard  for a better description of ${\cal P}_d$. \\

For every $(c,a):=(c_1,c_2,\cdot\cdot\cdot,c_{d-2},a)\in{\bf C}^{d-1}$ we denote by $P_{c,a}$ the polynomial of degree $d$ whose critical points are $(0,c_1,\cdot\cdot\cdot,c_{d-2})$
and such that $P_{c,a}(0)=a^d$. This polynomial is explicitely given by:

$$P_{c,a}:=\frac{1}{d} z^d+\sum_{2}^{d-1}\frac{(-1)^{d-j}}{j} \sigma_{d-j} (c) z^j +a^d$$

where $\sigma_i(c)$ is the symmetric polynomial of degree $i$ in $(c_1,\cdot\cdot\cdot,c_{d-2})$. 
For convenience we shall set $c_0:=0$.\\

We shall thus work within the holomorphic family 
$\big(P_{c,a}\big)_{(c,a)\in M}$ where the parameter space $M$ is simply ${\bf C}^{d-1}$. As explained in Milnor's paper \cite{Mi},
it is  convenient
to consider the projective compactification ${\bf P}^{d-1}$ of  ${\bf C}^{d-1}=M$ and see the sets $Per_n(w)$ as algebraic hypersurfaces of
${\bf P}^{d-1}$. We shall denote 
 the projective space at infinity $\{[c:a:0]\; ;(c,a)\in{\bf C}^{d-1}\setminus\{0\}\}$ by ${\bf P}_{\infty}$. \\

\subsection{The behaviour of the bifurcation locus at infinity}

We aim to show that  the bifurcation locus of the family $\{P_{c,a}\}_{ (c,a)\in{\bf C}^{d-1}}$ can only cluster on certain hypersurfaces of
${\bf P}_{\infty}$. The ideas here are essentially those used by Branner and Hubbard for proving the compactness of the connectedness locus (see \cite{BH} Chapter 1, section 3)
but we also borrow to the paper
(\cite{DF}) of Dujardin and Favre.\\

For every $0\le i\le d-2$, we will denote by $\alpha_i$ the homogeneous polynomial defined by 
$$\alpha_i (c,a):=P_{c,a}(c_i)=\frac{1}{d}c_i^d+\sum_{j=2}^{d-1}\frac{(-1)^{d-j}}{j} \sigma_{d-j} (c) c_i^j +a^d$$

and will consider the hypersurface $\Gamma_i$ of ${\bf P}_{\infty}$ defined by 
$$\Gamma_i:=\{[c:a:0]/\; \alpha_i(c,a)=0\}.$$

By a simple degree argument one sees that $P_{c,a}(0)=P_{c,a}(c_1)=\cdot\cdot\cdot=P_{c,a}(c_{d-2})=0$ implies that $c_1=\cdot\cdot\cdot=c_{d-2}=a=0$. This observation 
and Bezout's theorem lead to the following:

\begin{lem}\label{leminter}
The intersection $\Gamma_0\cap\Gamma_1\cap\cdot\cdot\cdot\cap\Gamma_{d-2}$ is empty and
$\Gamma_{i_1}\cap\cdot\cdot\cdot\cap\Gamma_{i_k}$ has codimension $k$ in ${\bf P}_{\infty}$ if $0\le i_1<\cdot\cdot <i_k \le d-2$.
\end{lem}

We shall denote by ${\cal P}_i$ the set of parameters $(c,a)$ for which the critical point $c_i$ of $P_{c,a}$ has a bounded forward orbit (recall that $c_0=0$).
The announced result can now be stated as follows.

\begin{theo}\label{controlinfty}
For every $0\le i\le d-2$, the cluster set of ${\cal P}_i$ in ${\bf P}_{\infty}$ is contained in $\Gamma_i$ and, in particular, the connectedness locus is compact in 
${\bf C}^{d-1}$.
\end{theo}

Since any cycle of attracting basins capture a critical orbit, the above theorem implies that the intersection of ${\bf P}_{\infty}$with an algebraic subset of the form $Per_{m_1}(\eta_1)\cap\cdot\cdot\cdot\cap Per_{m_k}(\eta_k)$ is contained
in some $\Gamma_{i_1}\cap\cdot\cdot\cdot\cap\Gamma_{i_k}$ as soon as the $m_i$ are mutually distinct and the $\vert\eta_i\vert$ strictly smaller than $1$.
Then, using Bezout's theorem, one gets the following:

\begin{cor}\label{interPer}
If $1\le k\le d-1$, $m_1<m_2<\cdot\cdot\cdot< m_k$ and $\sup_{1\le i\le k}\vert \eta_i\vert <1$ then 
$Per_{m_1}(\eta_1)\cap\cdot\cdot\cdot\cap Per_{m_k}(\eta_k)$ is an algebraic subset of codimension $k$ whose intersection with ${\bf C}^{d-1}$ is not empty. 
\end{cor}

The proof of theorem \ref{controlinfty} relies on estimates on the Green function and, more precisely, on the following result which is proved in 
the subsection $6.1$ of 
\cite{DF}.

\begin{prop} \label{estimgreen} Let  
$g_{c,a}(z):=\lim_nd^{-n}\ln^+ \vert P_{c,a}^n (z)\vert $ be the Green function of  $P_{c,a}$ and $G$ be 
 the function  defined on ${\bf C}^{d-1}$ by:
$G(c,a):=Max \{g_{c,a}(c_k);\;0\le k\le d-2\}$.
Let $\delta:= \frac{\sum_{k=0}^{d-2}c_k}{d-1}$. Then the following estimate occur:
\begin{itemize}
\item[1)] $Max \{g_{c,a}(z),G(c,a)\}\ge \ln \vert z-\delta\vert -\ln 4$
\item[2)] $G(c,a)=\ln^+ Max\{\vert a\vert,\vert c_k\vert\} +O(1).$
\end{itemize}
\end{prop}

\noindent {\it Proof of theorem \ref{controlinfty}.}  
Let $\Vert(c,a)\Vert:=Max\big(\vert a\vert,\vert c_k\vert\big)$. We simply have to check that  $\alpha_i\big(\frac{(c,a)}{\Vert (c,a)\Vert}\big)$
tends to $0$ when $\Vert(c,a)\Vert$ tends to $+\infty$ and $g_{c,a}(c_i)$ stays equal to $0$.
As $P_{c,a}(c_i)=\alpha_i(c,a)$ and $g_{c,a}(c_i)=0$, the estimates given by the proposition \ref{estimgreen} yield :
\begin{eqnarray*}
\ln^+\Vert (c,a)\Vert +O(1)=Max\big(dg_{c,a}(c_i),G(c,a)\big)=Max\big(g_{c,a}\circ P_{c,a} (c_i),G(c,a)\big)\ge \\
\ge\ln\frac{1}{4}\vert\alpha_i(c,a)-\delta\vert 
\end{eqnarray*}
 since $\alpha_i$ is $d$-homogeneous we then get for $\Vert (c,a)\Vert>1$:
$$(1-d)\ln \Vert (c,a)\Vert +O(1) \ge \ln\frac{1}{4}\vert\alpha_i\big(\frac{(c,a)}{\Vert (c,a)\Vert}\big)-\frac{\delta}{\Vert (c,a)\Vert^d}\vert$$
and the conclusion follows since $\frac{\delta}{\Vert (c,a)\Vert^d}$ tends to $0$ when $\Vert(c,a)\Vert$ tends to $+\infty$.\qed

\subsection{Proof of the main result}\label{proof}

We shall denote by $\la$ the parameter in ${\bf C}^{d-1}$ (i.e. $\la:=(c,a)$) and will then set $$\varphi_n(\la):=d^{-n}\ln \vert p_n(\la,w)\vert$$ where the polynomials $p_n(\la,w)$ are those given by the theorem \ref{theopoly}. 
We have to show that the sequence $(\varphi_n)_n$ is converging to $L$ in $L_{loc}^1$. 
When $\vert w\vert<1$, this has been shown to be true for any holomorphic family of rational maps (see the first assertion of theorem \ref{theoappL}), so we assume that $\vert w\vert=1$.\\

As it has been previously observed, the case $d=2$ is covered by the proposition \ref{propmodel}. To give a flavour of the proof when $d\ge 2$, we will first sketch it for $d=3$. \\

 {\bf Sketch of proof for degree three polynomials}. Let us first treat the problem on a curve $Per_{k_0}(\eta_0)$ for $\vert \eta_0\vert <1$. We will show that the sequence  $\varphi_n(\la)$ is uniformly converging to $L$ near
any stable (in $Per_{k_0}(\eta_0)$) parameter $\la_0$. To this purpose, one
desingularizes an irreducible component of $ Per_{k_0}(\eta_0)$ containing $\la_0$ and thus obtains a one-dimensional holomorphic family
$(P_{\pi(u)})_{u\in M}$. Keeping in mind that the elements of this family are degree 3 polynomials which do admit an attracting basin of period $k_0$ and using the fact that the 
connectedness locus in ${\bf C}^2$ is compact, one sees that the family $(P_{\pi(u)})_{u\in M}$ satisfies the assumptions of proposition \ref{propmodel}. The associated sequence $L_n(u,w)=\varphi_n(\pi(u))$
is therefore converging in $L^1_{loc}$ to $L$ and this convergence is locally uniform on stable components by proposition \ref{propelem}. \\
Let us now consider the problem on the full parameter space ${\bf C}^2$. Since the family  $\{P_{c,a}\}_{ (c,a)\in{\bf C}^{2}}$ contains hyperbolic parameters,
the first assertion of proposition \ref{prophyp} shows that the sequence $\varphi_n(\la)$ does not converge to $-\infty$.
According to the theorem \ref{compapsh}, it thus suffices to show that 
any $p.s.h$ limit value $\varphi$ of $\varphi_n(\la)$ in $L^1_{loc}({\bf C}^2)$ coincides with $L$.
Let us therefore assume that $\varphi_{n_k}$ tends to $\varphi$ in $L^1_{loc}({\bf C}^2)$.\\
We first show that $\varphi= L$ on any open subset of the type $A_{k_0}:=\cup_{\vert \eta\vert<1} Per_{k_0}(\eta)$.  According to the second assertion of proposition \ref{propelem}, it suffices to show that
$\varphi=L$ on any stable component $\Omega$ of $A_{k_0}$. By the third assertion of proposition \ref{propelem},
$\varphi_{n_k}$ is actually converging  pointwise to $\varphi$
on $\Omega$. As, by the previous step, $\varphi_n(\la)$ converges locally uniformly on the stable components of $Per_{k_0}(\eta)$, one thus obtains that $\varphi=L$ on $\Omega$.\\
According to the theorem \ref{controlinfty}, the set of non-hyperbolic parameters in ${\bf C}^2$ can only cluster on a finite subset of
${\bf P}_{\infty}$. We may therefore foliate ${\bf C}^2$ by parallel complex lines $(T_t)_{t\in{\bf C}}$ whose intersection with the  set of non-hyperbolic parameters is compact.
After taking a subsequence we may asume that $\varphi_{n_k}$ in converging to $\varphi$ in $L^1_{loc}(T_t)$ for almost every $t\in{\bf C}$.
To conclude it remains to see that $\varphi\vert_{T_t}\equiv L\vert_{T_t}$ for these $t$. For this, one uses lemma \ref{lemmodel}. The assumptions of the lemma are satisfied since, by construction, 
the unbounded stable component of $T_t$ is hyperbolic and the bifurcation locus in $T_t$ is accumulated by sets of the form $T_t\cap A_{k_0}$ where, as we have previously shown, $\varphi=L$.\qed\\

{\bf Proof of the theorem \ref{theoequipol}}.
For $1\le q\le d-2$, the notation $W_q$ will refer to any irreducible component of  a $q$-codimensional analytic subspace of  ${\bf C}^{d-1}$  of the form $Per_{n_1}(\eta_1)\cap\cdot\cdot\cdot\cap Per_{n_q}(\eta_q)$
where $(\eta_1,\cdot\cdot\cdot,\eta_q)\in \Delta^q$ and the integers $n_j\ge 2$ are mutually distinct (by corollary \ref{interPer} such sets do exist). Let us stress that if $\la\in W_q$ then the polynomial $P_{\la}$ admits $q$ distinct attracting basins
besides the basin at infinity.
Analogously we shall set $W_0:={\bf C}^{d-1}$. By $W_q^{reg}$ we shall denote the regular part of $W_q$. The proof will consist in showing  by decreasing induction on $ 0\le q\le d-2$ that 
$$(*_q):\;\;\textrm{the sequence}\; \varphi_n \vert_{W_q}\;\textrm{tends to}\; L \;\textrm{in}\; L^1_{loc}(W_q^{reg})\;\textrm{for any}\; W_q.$$

Let us first establish $(*_{d-2})$. The analytic set $W_{d-2}$ is a curve in ${\bf C}^{d-1}$. Desingularizing we get
a proper holomorphic map $\pi:M\to W_{d-2}$ where $M$ is a Riemann surface. We claim that the one-dimensional holomorphic family $(P_{\pi(u)})_{u\in M}$ satisfies the assumptions of proposition \ref{propmodel}.
To see this we observe that there exists at most one critical point of the polynomial $P_{\pi(u)}$ whose orbit is not captured by one of the $d-2$ distinct attracting basins of $P_{\pi(u)}$. Let us denote
by $c(u)$ this critical point. Assume that $u_0$ lies in the bifurcation locus of $(P_{\pi(u)})_{u\in M}$. Since all critical points of $P_{\pi(u)}$, except a priori $c(u)$, stay in some attracting basin for $u$ close to $u_0$, the orbit of $c(u)$ cannot be uniformly bounded on a small neighbourhood of $u_0$. 
This implies that $c(u)$ must belong to the basin of infinity for a convenient small perturbation of $u_0$ and shows that $P_{\pi(u_0)}$ becomes hyperbolic after a convenient small pertubation . In other words, the bifurcation locus of
$(P_{\pi(u)})_{u\in M}$ is accumulated by hyperbolic parameters.\\
The above argument also shows that
if $P_{\pi(u)}$ is non-hyperbolic, then $c(u)$ cannot belong to the basin at infinity and therefore $P_{\pi(u)}$ belongs to the connectedness locus.
Using the compactness of the connectedness locus and the properness of the map $\pi$, one sees that the set of non-hyperbolic parameters of $M$ is compact.\\
By proposition \ref{propmodel}, $\varphi_n(\pi(u))$ is converging in 
$L^1_{loc}(M)$ to $L\circ \pi$. By the third assertion of the proposition \ref{propelem}, the convergence is actually pointwise on the stable components of $M$ and thus
$\varphi_n$ converges pointwise to $L$ on the stable set of $W_q^{reg}$. By the second assertion of 
proposition \ref{propelem}, 
 $\varphi_n \vert_{W_q}\;\textrm{  converges to}\; L \;\textrm{in}\; L^1_{loc}(W_q^{reg})$.
We have proved $(*_{d-2})$.\\

Assuming now that $(*_{q+1})$ is satisfied, we shall prove that $(*_q)$ is true. Let us fix an irreducible $q$-codimensional analytic set 
$W_q\subset Per_{n_1}(\eta_1)\cap\cdot\cdot\cdot\cap Per_{n_q}(\eta_q)$.\\ 
One easily deduce from corollary \ref{interPer} that $W_q$ contains hyperbolic parameters and this fact
preserves
$\varphi_{n}\vert_{W_q^{reg}}$ to converge to $-\infty$ (see proposition \ref{prophyp}).
According to theorem \ref{compapsh}, we thus have to show that for
any subsequence $\varphi_{n_k}\vert_{W_q^{reg}}$ converging
to some $p.s.h$ function $\varphi$ in $L^1_{loc}(W_q^{reg})$ one actually has $\varphi= L\vert_{W_q^{reg}}$.\\

We shall use the two following facts which will be proved later.\\

{\bf Fact 1} {\it  Let $A_m$ be an open subset of ${\bf C}^{d-1}$ defined by $A_m:=\cup_{\vert \eta\vert<1} Per_m(\eta)$ where $m > max(n_1,\cdot\cdot\cdot,n_q)$. 
If $W_q^{reg}\cap A_m$ is not empty, then $\varphi=L$ on $W_q^{reg}\cap A_m$.}\\

{\bf Fact 2} {\it  There exists a foliation $\cup_{t\in A} T_t$ of ${\bf C}^{d-1}$ by $(q+1)$-dimensional parallel affine subspaces such that, for almost every $t\in A$, the slices $T_t \cap W_q$ are curves
on which the set of non-hyperbolic parameters is relatively compact.}  \\

Let us consider the curves $T_t \cap W_q$ which are given by Fact 2. By Slutsky lemma, $\varphi_{n_k}$ is converging to $\varphi$ in $L_{loc}^1$ on almost all these curves and it
thus remains to show that $\varphi=L$ on them. To this purpose we consider an irreducible component $\Gamma$ of $T_t\cap W_q$ and desingularize it. This yields a proper holomorphic map 
$\pi :M\to \Gamma$ where $M$ is a Riemann surface. We shall reach the conclusion by
applying lemma \ref{lemmodel} to the family $\big( P_{\pi(u)}\big)_{u\in M}$.\\
By the properness of $\pi$ and Fact 2, the set of non-hyperbolic parameters in $M$ is compact and therefore any non relatively compact stable component in $M$ is hyperbolic. Then, by the first assertion of proposition \ref{prophyp}, $\varphi\circ\pi=L\circ \pi$ on such components.\\
Using Fact 1, we shall now prove that the bifurcation locus of $\big( P_{\pi(u)}\big)_{u\in M}$ is accumulated by parameters where $\varphi\circ\pi =L\circ \pi$. Let $u_0$ be a point in the bifurcation locus,
 we may assume that $\pi$ is locally biholomorphic at $u_0$ and it thus suffices to accumulate 
$\pi (u_0)$ by points where $\varphi =L$. 
As it is well known, $u_0$ is accumulated by parameters 
$u_k$ such that 
$P_{\pi (u_k)}\in Per_{m_k}(0)$ 
and $m_k\to +\infty$ (this follows also from the general fact that $T_{bif}=\lim_m d^{-m} [Per_m(0)]$).  This implies that
$\pi (u_0)$ is accumulated by open sets of the form $W_q\cap A_{m_k}$. It then follows from Fact 1 that $\pi (u_0)$ is accumulated  by points $\la_k$ for which $\varphi (\la_k)=L(\la_k)$. This ends the proof.\\

Let us finally establish the Facts. \\
{\it Fact 1.} Let $\Omega$ be a stable component of $W_q^{reg}\cap A_m$. According to the first and third assertions of proposition \ref{propelem}, the sequence
$\varphi_{n_k}-L$ is locally uniformly converging to the pluriharmonic negative function $\varphi -L$
on $\Omega$ (as previously observed, $W_q$ contains hyperbolic parameters and therefore has no persistent neutral cycles). 
For all but a finite number of $\eta\in \Delta$ the analytic set  $W_q\cap Per_m(\eta)$ is of codimension $q+1$ (otherwise $W_q$ would be contained in infinitely many hypersurfaces $Per_m(\eta)$  and $P_{\la}$ would have an infinite number of attracting basins when $\la\in W_q$). 
 Let us thus pick $\eta_0\in\Delta$ and $\la_0\in \Omega\cap Per_m(\eta_0)$  such that  $W_q\cap Per_m(\eta_0)$ has codimension $q+1$ and is regular at $\la_0$. Let us denote by $W_{q+1}$
 the irreducible component of $W_q\cap Per_m(\eta_0)$ to which belongs $\la_0$.
 Then, by construction, $\la_0$ belongs to some stable component $\omega$ of $W_{q+1}^{reg}$. Combining the induction assumption $(*_{q+1})$, 
with the third assertion of proposition \ref{propelem} on sees that $\varphi-L=0$ on $\omega$. In particular $\varphi(\la_0) -L(\la_0)=0$ and, by the maximum principle,
$\varphi -L=0$ on $\Omega$.\\
It now follows from 
 the second assertion of proposition \ref{propelem} that $\varphi=L$ on $W_q^{reg}\cap A_m$. The Fact 1 is proved.\\
\\
{\it Fact 2.} Let $\widetilde{W_q}$ be the algebraic subset of ${\bf P}^{d-1}$ such that $\widetilde{W_q}\cap {\bf C}^{d-1}=W_q$. When $q>0$ and $\la \in W_q$ then $P_{\la}$ has $q$ distinct attracting basins and, therefore, at least $q$ of its critical points have a bounded orbit. According to theorem \ref{controlinfty} we thus have
$$\widetilde{W_q}\cap {\bf P}_{\infty} \subset \bigcup_{0\le i_1<\cdot\cdot\cdot<i_q\le d-2}\Gamma_{i_1}\cap \cdot\cdot\cdot \cap\Gamma_{i_q}$$
and, moreover, $\bigcup_{0\le i_1<\cdot\cdot\cdot<i_{q+1}\le d-2}\Gamma_{i_1}\cap \cdot\cdot\cdot \cap\Gamma_{i_{q+1}}$
is a $(d-3-q)$-dimensional algebraic subspace of  ${\bf P}_{\infty}$.
Thus,  as it is classical (see \cite{Chi} 7.3),
we may pick a $q$-dimensional complex plane $C_{\infty}$ in ${\bf P}_{\infty}$ (a point when $q=0$) such that 
$$C_{\infty}\cap\big(\bigcup_{0\le i_1<\cdot\cdot\cdot<i_{q+1}\le d-2}\Gamma_{i_1}\cap \cdot\cdot\cdot \cap\Gamma_{i_{q+1}}\big) =\emptyset.$$
We now slice ${\bf C}^{d-1}$ by $(q+1)$-dimensional parallel affine subspace $T_t$ which cluster on $C_{\infty}$ in ${\bf P}^{d-1}$ and
write ${\bf C}^{d-1}=\cup_{t\in A} T_t$ where $A$ is a $(d-q-2)$-dimensional complex plane
which is transverse to the foliation.\\
If $\la\in W_q$ then at least $q$ of the critical points of $P_\la$ belong to some attracting basin. This implies that
 the set of non-hyperbolic parameters in $W_q\cap T_t$  may only cluster on the intersection of $C_{\infty}$ with 
$\cup_{0\le i_1<\cdot\cdot\cdot<i_{q+1}\le d-2}\Gamma_{i_1}\cap \cdot\cdot\cdot \cap\Gamma_{i_{q+1}}$.
 The choice of $C_{\infty}$ 
 guarantees therefore that, for all $t\in A$, the set of non-hyperbolic parameters in $W_q\cap T_t$ is compact.\\
It remains to show that, for almost all $t\in A$, the analytic set $W_q\cap T_t$ is a curve. Let us, to this purpose, denote by $\sigma: W_q\to A$ the canonical projection from $W_q$ onto $A$.
The fibers of $\sigma$ are the analytic sets $W_q\cap T_t$ whose dimensions are at least equal to $(d-1)-dim A -q=1$. Then, the set of points $a\in A$ for which the fiber $\sigma^{-1} \{a\}$ 
is of dimension strictly greater than $1$ is contained in a countable union of analytic subsets of $A$ whose dimensions are smaller than $dim W_q -2= (d-1) -q -2 =dim A -1$ 
(see \cite{Chi}, 3.8) and is therefore
Lebesgue negligeable. In other words $W_q\cap T_t$ is a curve for almost all $t$ and Fact 2 is proved.\qed

\bibliographystyle{amsalpha}

\end{document}